\documentclass[11pt]{amsart}
\usepackage{times,mathptmx,graphics,amssymb,amsthm}

\newtheorem{theorem}{Theorem}[section]
\newtheorem{theorema}{Theorem}

\newtheorem{proposition}[theorem]{Proposition}
\newtheorem{lemma}[theorem]{Lemma}

\newtheorem{conjecture}{Conjecture}

\newtheorem{corollary}[theorem]{Corollary}

\theoremstyle{remark}
\newtheorem*{remark}{Remark}

\newcommand{\fig}[2]{
    \begin{figure}[ht]
    \begin{center}
    \includegraphics{#1.eps}
    \end{center}
    \caption{#2}\label{#1}
    \end{figure}}

\date{November 24, 1995}

\title{Highly saturated packings and reduced coverings}
\author{Gabor Fejes T\'oth, Greg Kuperberg, and W{\l}odzimierz Kuperberg}

\begin{document}

\maketitle

\begin{abstract}
We introduce and study certain notions which might serve as substitutes for
maximum density packings and minimum density coverings. A body is a compact
connected set which is the closure of its interior.  A packing $\mathcal{P}$ with
congruent replicas of a body $K$ is $n$-saturated if no $n-1$ members of it
can be replaced with $n$ replicas of $K$, and it is completely saturated if
it is $n$-saturated for each $n\ge 1$.  Similarly, a covering $\mathcal{C}$ with
congruent replicas of a body $K$ is $n$-reduced if no $n$ members of it can
be replaced by $n-1$ replicas of $K$ without uncovering a portion of the
space, and it is completely reduced if it is $n$-reduced for each $n\ge 1$.
We prove that every body $K$ in $d$-dimensional Euclidean or hyperbolic
space admits both an $n$-saturated packing and an $n$-reduced covering with
replicas of $K$.
Under some assumptions on $K\subset \mathbb{E}^d$
(somewhat weaker than convexity), we prove the existence of completely
saturated packings and completely reduced coverings, but in general, the
problem of existence of completely saturated packings and completely reduced
coverings remains unsolved.
Also, we investigate some problems related to the the
densities of $n$-saturated packings and $n$-reduced coverings. Among other
things,
we prove that there exists an upper bound for the density of a $d+2$-reduced
covering of $\mathbb{E}^d$ with congruent balls, and we produce some density
bounds for the $n$-saturated packings and $n$-reduced coverings of the plane
with congruent circles.
\end{abstract}

\section{Introduction and preliminaries}

Two of the basic problems in the theory of packing and covering are to
determine the most efficient packing and covering with replicas of (meaning
sets congruent to) a given set $K$ in some metric space.  Recall that a
packing is a family of sets whose interiors are mutually disjoint, and that a
covering is a family of sets whose union is the whole space.
By a space we mean either $d$-dimensional Euclidean space $\mathbb{E}^d$
or $d$-dimensional hyperbolic space $\mathbb{H}^d$, although the definitions that
follow are sometimes more general.  We shall consider packings and coverings
with replicas of a nonempty compact connected set which is the closure of its
interior, a {\it body}, for short.

The usual measure of the efficiency of an arrangement in Euclidean space is
density. Roughly speaking, the density of an arrangement is the total volume
of the members of the arrangement divided by the volume of the whole space. 
Rigorously, density can be defined by an appropriate limit \cite{FTK}. The
maximum density of a packing of the space with replicas of a
(measurable) set $K$ is denoted by $\delta (K)$ and is called the {\it
packing density of} $K$. The minimum density of a covering with replicas
of $K$ is denoted by $\vartheta (K)$ and is called the {\it covering
density of} $K$.  It is known that each of the maximum and the minimum
density is attained \cite{G}.

There are some disadvantages of using density for measuring the efficiency of
an arrangement. In the first place, optimum density is a global notion and it
does not imply the local efficiency of an arrangement. Secondly, the notion
of density cannot be extended in full generality to hyperbolic geometry
\cite{Bor1}, \cite{FTK}. In what follows we introduce and study certain
notions which might serve as substitutes for maximum density packings and
minimum density coverings.

Let $K$ be a body and let $\mathcal P$ be a packing of space with replicas of $K$.
$\mathcal P$ is said to be {\it saturated} if it cannot be augmented with any
additional replica of $K$ without overlapping with a member of $\mathcal P$. More
generally, $\mathcal P$ is $n$-{\it saturated} if no $n-1$ members of it can be
replaced with $n$ replicas of $K$. A packing is {\it completely saturated} if
it is $n$-saturated for every $n\ge 1$.

Note that Fejes T\'oth and Heppes \cite{FTH} define the term ``$n$-saturated''
differently, but we hope that our definition causes no confusion.

A covering $\mathcal C$ of $E^d$ with replicas of $K$ is {\it reduced} if no
proper sub-family of $\mathcal C$ is a covering. Similarly,
we say that $\mathcal C$ is $n$-{\it reduced} if no
$n$ members of it can be replaced by $n-1$ replicas of $K$ without uncovering
a portion of the space. A covering is {\it completely reduced} if it is
$n$-reduced for every $n\ge 1$.
 
\begin{conjecture} Every body $K$ in $\mathbb{E}^d$ (resp. in $\mathbb{H}^d$)
admits a completely saturated packing and a completely reduced
covering with replicas of $K$.  \label{ccomplete}
\end{conjecture}

This conjecture is supported by the following results:
 
\begin{theorem} Every convex body in $\mathbb{E}^d$ admits a completely saturated
packing and a completely reduced covering of $\mathbb{E}^d$ with replicas of the
body. \label{thcomplete}
\end{theorem}

\begin{theorem} Every body $K$ in $\mathbb{E}^d$ (resp. in $\mathbb{H}^d$)
admits both an $n$-saturated packing and an $n$-reduced covering with
replicas of $K$. \label{thn}
\end{theorem}

Section~\ref{scomplete} presents a proof of Theorem~\ref{thcomplete}.
We note there that the theorem holds for bodies satisfying the strict
nested similarity property, a condition weaker than convexity.
Theorem~\ref{thn} is proved in Section~\ref{slattice} (for the Euclidean
case) and Section~\ref{shyperbolic} (for the hyperbolic case), each
as a corollary of a more general statement.  The hyperbolic case involves
some elements of the theory of hyperbolic manifolds, which we review for
the benefit of the unfamiliar reader.

\subsection{Acknowledgements}

The authors acknowledge, with gratitude, that during the preparation
of this paper, the research of G. Fejes T\'oth has been supported by the
Hungarian Foundation for Scientific Research (OTKA), grants no.~1907 and
no.~14218, and that of W. Kuperberg by the National Science Foundation,
grant no.~DMS-9403515.

\section{Complete saturation and reduction}
\label{scomplete}

In this section, we give a proof of Theorem~\ref{thcomplete}.
We precede the proof with some definitions and two lemmas.
Throughout the argument, $K$ is a given convex body in $\mathbb{E}^d$, and, as
before, $V(A)$ denotes the volume of $A$. We use the Hausdorff distance
between closed sets to measure the distance between a pair of (finite)
packings or coverings, extending the Hausdorff distance function to the
space of finite (unordered) collections of compact sets in the natural way.

Let $c$ be a ``center'' point in the interior of $K$, say the center of
gravity of $K$. Let $B(r,p)$ be the sphere of radius $r$ centered at
$p$.  A packing with replicas of $K$  is {\it completely saturated in
$B(r,p)$} if no $n$ replicas contained in $B(r,p)$ can be replaced by
$n+1$ replicas contained in $B(r,p)$, for every integer $n$ (the
replicas not contained in $B(r,p)$ are not to be moved in this
process).  A packing is {\it unsaturated in $B(r,p)$} for short if it
fails to be completely saturated in $B(r,p)$.  An arbitrary packing of
$\mathbb{E}^d$ can be altered within $B(r,p)$ so as to result in a
packing completely saturated in the ball: delete all replicas of $K$
contained in the ball and replace them with the maximum number of
replicas of $K$ that will fit in the ball without overlapping with each
other or with any of the replicas that partially invade the ball.

A {\it homothetic thinning} $T_h(\mathcal P)$ of a packing $\mathcal P$ with replicas
of $K$ by a factor of $h>1$ is a new packing such that each center $c\in
\mathbb{E}^d$ maps to $hc$, but such that the replicas of $K$ are translated
without expansion. The analogous concepts for coverings ({\it completely
reduced in $B(r,p)$}, {\it unreduced in $B(r,p)$}, {\it homothetic
thickening} by a factor of $h<1$) are defined similarly. The proof of
Theorem~\ref{thcomplete} is given only for packings, since the proof for
coverings is the same except for one modification which is mentioned
afterwards.

\begin{lemma} For every $r$ and $\varepsilon>0$, there exists a $\delta>0$
such that if a $\mathcal P$ is less than $\delta$ away (in Hausdorff distance)
from a packing which is unsaturated in $B(r,0)$, then
$T_{1+\varepsilon}(\mathcal P)$ is unsaturated in $B((1+\varepsilon)r,0)$.
\label{lhausdorff} \end{lemma}

The proof is left to the reader.

\begin{lemma} Let $r>0$ and $\eta>0$.  Then there exists an $s_0$ and a
$\delta>0$ such that for every $s > s_0$, a packing $\mathcal P$ of replicas of
$K$ which is densest relative to $B(2s+r,0)$ has the following property: If
$p$ is chosen at random in $B(s,0)$, $\mathcal P$ is at least $\delta$ away from
unsaturated in $B(r,p)$ with probability at least $1-\eta$.
\label{lprobability} \end{lemma}

\begin{proof} Informally, if $\varepsilon$ is sufficiently small and $s$ is
sufficiently large, then if $\mathcal P$ is expanded by $1+\varepsilon$, the
loss of density from replicas of $K$ sliding over the edge of $B(2sR+r,0)$ is
outweighed by the gain in finding a $\eta$ proportion of $B(r,p)$'s inside
that are unsaturated and re-saturating the packing in a disjoint collection
of these smaller balls. Then $\delta$ can be chosen based on $\varepsilon$
and Lemma~\ref{lhausdorff}.  A more precise argument follows.

Temporarily fix $\delta > 0$ and $s>0$, and suppose that, to the contrary,
the set $X$ of points $p \in B(s,0)$ such that $\mathcal P$ is less than
$\delta$ away from unsaturated in $B(r,p)$ has measure at least $\eta
V(B(s,0))$.  We will arrive at a contradiction for $\delta$ sufficiently
small and $s$ sufficiently large.

Since $X$ has measure $\eta V(B(s,0))$, it cannot be covered by fewer than
$$\eta {V(B(s,0)) \over V(B(2r,0)) } = {\eta \over (2r)^d}s^d$$ balls of
radius $2r$.  It follows that there exists a packing $\{B(r,p_i)\}_{1 \le
i \le k}$ of $k$ balls of radius $r$ entirely within $B(s+r,0)$ such that
the restriction of $\mathcal P$ to each ball is less than $\delta$ away from
an unsaturated packing, where $$k > \eta {V(B(s,0)) \over V(B(2r,0))} >
cs^d$$ for some constant $c$ depending only on $r$ and $\eta$. In other
words, $${k V(K) \over V(B(s,0))} > {cV(K) \over V(B(1,0))}.$$

Let $$C={cV(K)\over V(B(1,0))} < 1,$$ let
$$\varepsilon=\min\left((1-C)^{-{1\over d}}-1,\frac12\right),$$
and let $\delta$ be given by Lemma~\ref{lhausdorff}.

Observe that the difference between the density of $\mathcal P$ and that of the
homothetic thinning $T_{1+\varepsilon}(\mathcal P)$ (both relative to $B(s,0)$)
is at most $1-(1-\varepsilon)^{-d}$. On the other hand, since
$T_{1+\varepsilon}(\mathcal P)$ is unsaturated in each of the balls
$B((1+\varepsilon)r,(1+\varepsilon)p_i)$, the density of
$T_{1+\varepsilon}(\mathcal P)$ relative to $B(s,0)$ can be increased through
saturation by an amount greater than
$$k(s)V(K)/V(B(s,0)>cV(K)/V(B(1,0))=C.$$ (Note that by our choice of
$\varepsilon$, each $B((1+\varepsilon)r,(1+\varepsilon)p_i)$ is contained
in $B(2s+r,0)$.) By our choice of $\varepsilon$, we have
$1-(1-\varepsilon)^{-d}<C$. Thus, the thinned and then re-saturated packing
is denser in $B(2s+r,0)$ than the original packing $\mathcal P$, which is a
contradiction. \end{proof}

\begin{proof}[Proof of theorem] Let $n>0$ be an integer and for each $1 \le k \le n$, let
$\delta_k$ and $s_k$ be given by the second lemma with $r = k$ and $\eta
=3^{-k}$.  Let $s$ be the supremum of the $s_k$'s, and let $\mathcal P$ be a
packing which is densest relative to $B(2s+n,0)$.  Then for at least
$$1-(\frac13 + \frac19 + \ldots + \frac1{3^n}) > 1/2$$ of $p \in B(s,0)$,
$\mathcal P$ is $\delta_k$ away from unsaturated in all $B(k,p)$'s
simultaneously. After translation by $-p$, $\mathcal P$ becomes a packing ${\mathcal
P}_n$ which is simultaneously $\delta_k$ away from unsaturated in each
$B(k,0)$. The sequence ${\mathcal P}_n$ has a subsequence which converges in the
Hausdorff topology to a limit $\widetilde{\mathcal P}$. The packing
$\widetilde{\mathcal P}$ is $\delta_n$ away from unsaturated in $B(n,0)$ for every
$n > 0$, and it is therefore completely saturated in $\mathbb{E}^d$.
\end{proof}

{\sc Remark 1.} The argument for the analogous theorem for coverings
requires two minor modifications in the formulation and proof of
Lemma~\ref{lprobability}. Firstly, instead of using a covering which
extremizes density in $B(2s+r,0)$, we use a minimum cardinality arrangement
of replicas of $K$ that covers $B(2s+r,0)$. In particular, an optimal
covering of $B(2s+r,0)$ in this sense has no replicas disjoint from
$B(2s+r,0)$. Secondly, a homothetic thickening of a covering of $B(2s+r,0)$
by a factor of $1-\varepsilon$ is not in general again a covering. To repair
it, we identify a cube $\kappa \subset K$ and we cover the annular region
$B(2s+r,0) - B((1-\varepsilon)(2s+r),0)$ by non-overlapping translates of
$\kappa$. The resulting gain in number of replicas due to homothetic
thickening is comparable to the loss of density due to homothetic thinning.
\bigskip

{\sc Remark 2.} The hypothesis that $K$ is convex can be weakened
somewhat without any changes in the proof. It suffices that $K$ have
the {\it strict nested similarity property}, which requires that, for
every positive number $h<1$, the interior of $K$ contains a replica of
$hK$. For example, if $K$ is {\it strictly starlike}, i.e., $K$
contains an interior point such that every ray emanating from it meets
the boundary of $K$ at a single point, then $K$ has the strict nested
similarity property. Another example of a body with the strict nested
similarity property is an $\varepsilon$-neighborhood of a logarithmic
spiral, as shown in Figure~\ref{figure4}.

\fig{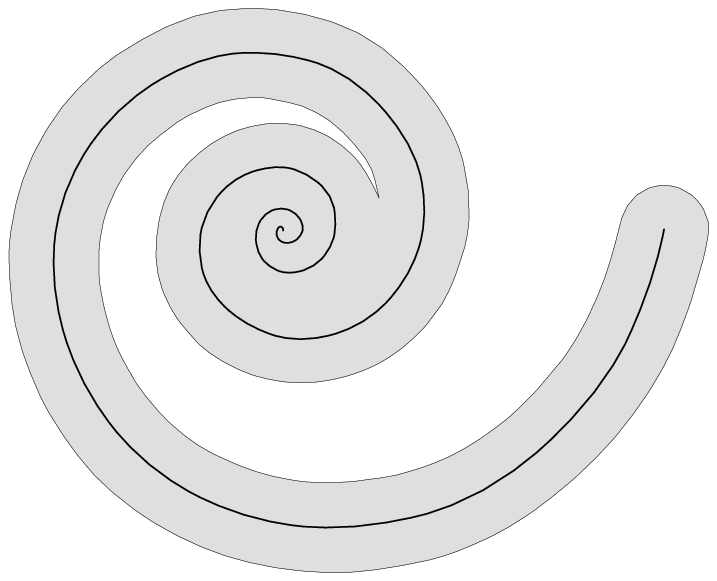}{A body with the strict nested similarity property}

\section{Lattices of isometries}
\label{slattice}

We recall a construction from the topological theory of covering spaces, to be
used in this section and in the next one, which we apply to produce
$n$-saturated packings and $n$-reduced coverings in various spaces. For the
basic notions and facts related to that theory we refer the reader to
\cite[Ch. 1-2]{S}.

Let $M$ be a locally compact, connected metric space (for our purposes it is
sufficient to assume that $M$ is a non-compact Riemannian manifold) with the
metric $\varrho$. A group $G$ of self-isometries of $M$ is a
{\it fixed-point-free, uniform lattice of isometries}, or
{\it lattice of isometries} or {\it lattice} for short, if $G$ satisfies the
following two conditions:

\begin{description}
\item[({\it i})] there is a number $\gamma>0$ such that for every $x\in
M$ and every $g\in G$ other than the identity,
$\varrho\left(x,g(x)\right)\ge\gamma$;

\item[({\it ii})] the quotient space $B=M/G$ (whose points are the
orbits $\{g(x): g\in G\}$ of points of $M$ under $G$, furnished with
the quotient topology) is compact.
\end{description}

Every lattice $G$ of isometries of $M$ has the crucial property that the
quotient map $p$ from $M$ to $B$, assigning to each point of $M$ its orbit,
is a covering projection (see \cite[p. 88]{S}).  We shall refer to $M$ as the
covering space and to $B$ as the base space. Also, the covering is regular,
for $G$ acts transitively on each point-inverse. The base space
is metrizable: a specific metric for $B$ is defined by
$\varrho(x,y)=\inf\{\varrho\left(\tilde x, \tilde y\right): p(\tilde x)=x,
p(\tilde y)=y \}$. Under this metric, the covering projection $p$ is a local
isometry:  the restriction of $p$ to any set of
diameter smaller than $\gamma$ is an isometry.

Conversely, given a regular covering $p:M\to B$, where $B$ is a compact
connected manifold, $M$ is connected and endowed with the metric lifted from
$B$, then the group of covering transformations of $M$ is a lattice of
isometries. Moreover, if $M$ is simply connected, then the lattice of
isometries of $M$ is isomorphic to the fundamental group $\pi_1(B)$ of the
base space $B$ (see \cite[Sec. 2.6]{S}).

If $M$ is $d$-dimensional Euclidean space, then a lattice of isometries which
consists of translations reduces to the classical concept of a lattice of
vectors, and the resulting base space is a $d$-dimensional torus. However,
there are other lattices of isometries even in $\mathbb{E}^2$. For instance,
the group of isometries of the Cartesian plane generated by a translation in
the $x$ direction and a glide-reflection in the $y$ direction is a lattice of
isometries, but the resulting base space is a Klein bottle and not a torus.

Define the {\it girth} of a lattice $G$ as the infimum of the distances
$\varrho\left(x,g(x)\right)$ over all non-identity elements $g\in G$. If $M$
is a simply-connected Riemannian manifold (such as $\mathbb{H}^d$), then the same
number is the infimum of the lengths of all non-trivial loops in the base
space, hence the name ``girth''.  For many manifolds, including
Euclidean and hyperbolic manifolds, the girth is twice the {\em injectivity
radius} of the quotient manifold, which is defined as the largest
$r$ such that no metric ball of radius $r$ overlaps itself.

Observe that if $S$ is a subset of $M$ whose diameter is smaller than
the girth of $G$, then the image $p(S)$ is a replica of $S$ in $B$.
Also, the set $p^{-1}(p(S))$ is the union of a discrete collection of
mutually disjoint replicas of $S$, namely it is the orbit of $S$ under
$G$. We call this discrete collection of replicas of $S$ a {\it
lifting} of $p(S)$ (in $M$).

\begin{theorem} Let $M$ be a locally compact connected metric space and let
$K$ be a body in $M$. If for every $c>0$, $M$ admits a lattice of girth
greater than $c$, then there exist an $n$-saturated packing of $M$  and an
$n$-reduced covering of $M$ with replicas of $K$, for every positive
integer $n$.  \label{thifgirth}
\end{theorem}
\begin{proof} We restrict our attention to packings only, since the case of
coverings is completely analogous.

Let $G$ be a lattice of isometries of $M$ of girth greater than $2n+1$ times
the diameter of $K$, let $B$ be the base space associated with the lattice,
and let $p:M\to B$ be the covering projection. In the base space $B$, arrange
a packing with a maximum number of bodies of the form $p(K')$, where $K'$ is
a replica of $K$ in $M$. The maximum is finite because $B$ is compact and $K$
has a non-empty interior. Let $\mathcal P$ be the lifting (in $M$) of this
packing. Obviously, $\mathcal P$ is a packing of $M$ with replicas of $K$. We
assert that $\mathcal P$ is $n$-saturated.

Suppose the contrary, and let $m\le n$ be the smallest positive integer such
that $\mathcal P$ is not $m$-saturated. Therefore there are $m$ members of
$\mathcal P$, say $K_1, K_2,\ldots , K_m$ which can be replaced by $m+1$ other
replicas of $K$, say $L_1, L_2, \ldots$, $L_{m+1}$, and  $m$ is the smallest
integer with this property. By the ``pigeonhole principle,'' the set
$S=\left(\bigcup K_i\right)\cup \left(\bigcup L_j\right)$ is connected.
Therefore the diameter of $S$ is smaller than the sum of the diameters of the
$K_i$'s and the $L_j$'s, thus smaller than the girth of $G$. If we replace in
$B$ the sets $p(K_i)\ (i=1,2,\ldots , m)$ with the sets $p(L_j)\
(j=1,2,\ldots , m+1)$, we exceed the maximum number defined above. This is a
contradiction.
\end{proof}

The Euclidean case of Theorem~\ref{thn} is an immediate corollary of
Theorem~\ref{thifgirth}.
Also, the remaining part of Theorem~\ref{thn} is now reduced to the problem of
existence of lattices of arbitrarily large girth in $d$-dimensional
hyperbolic space. This problem is addressed in the next section.

\section{Hyperbolic lattices of large girth}
\label{shyperbolic}

A lattice of isometries of hyperbolic space $\mathbb{H}^d$ will be called
a ($d$-dimensional) {\it hyperbolic lattice} for short.  The aim of
this section is to prove the following:

\begin{theorem}
For every $c$, there exists a $d$-dimensional hyperbolic lattice of
girth greater than $c$.  \label{thgirth}
\end{theorem}

Although this fact and the methods used for proving it have been known for a
long time, we could not find a suitable reference and we include a proof for
completeness.

We begin with some algebraic preliminaries. A group $G$ is {\it residually
finite} if for every $g\in G$ other than the identity $e$, there exists a
normal subgroup $N$ of finite index which does not contain $g$. Equivalently,
$G$ is residually finite if for every $g\in G,\ g\neq e$, there is a
homomorphism $\varphi$ from $G$ to a finite group such that $\varphi(g)$ is
not the identity.

Since the intersection of two normal, finite-index subgroups of $G$ is again a
normal subgroup of finite index, we get immediately:

\begin{proposition} If $G$ is residually finite, then for any finite set
$F\subset G$ not containing the identity there exists a normal subgroup
$N\subset G$ of finite index which does not intersect~$F$.
\label{pindex}
\end{proposition}

The group of non-singular $n\times n$ matrices with real coefficients is
denoted by $GL(n,\mathbb{R})$ and $I(\mathbb{H}^d)$ denotes the group of isometries
of $\mathbb{H}^d$. Since $\mathbb{H}^d$ can be modelled as one sheet of a two-sheeted
hyperboloid in $\mathbb{R}^{d+1}$, where the isometries of $\mathbb{H}^d$ are
those linear transformations of $\mathbb{R}^{d+1}$ which preserve the
sheet (see \cite[Sec. A2]{BP}), the group $I(\mathbb{H}^d)$ is isomorphic to a
subgroup of $GL(d+1,\mathbb{R})$.

The following lemma is a direct consequence of a theorem of Mal'cev
(\cite[Th.~VII]{M}). We include a version of Mal'cev's proof.

\begin{lemma} Every finitely generated subgroup $G$ of
$GL(n,\mathbb{R})$ is residually finite.  \label{lfinite}
\end{lemma}
\begin{proof} Let $g$ be a non-identity element of $G$. The aim of
the proof is to construct a homomorphism from $G$ to a finite group such that
the image of $g$ is also not the identity. This is accomplished by the
composition
of three homomorphisms.  The first one, $\alpha$, sends $G$ into the group of
algebraic matrices (matrices whose entries are algebraic numbers) of the same
size as the matrices in $G$; the second one, $\beta$, is a map to 
a group of (larger) rational matrices; and the third one, $\gamma$, is
a map to a group of matrices over a finite field $\mathbb{Z}/p$.

Let $g_1,g_2,\ldots ,g_k$ be a set of generators for $G$.
If a set of algebraic equations in finitely many variables has a real
solution, then it has (possibly complex) algebraic solutions arbitrarily
close to
the real solution. The defining relations between the $g_i$'s impose some
constraints on the entries of these matrices. Since the constraints are
algebraic, algebraic matrices can be found that satisfy the same relations as
the $g_i$'s do, and are arbitrarily close to them.  We pick algebraic
matrices for the images under $\alpha$ of the generators that are close
enough to the original real matrices so that $g$'s image $\alpha(g)$ is not
the identity.

The coefficients of all of the $\alpha(g_i)$'s are algebraic numbers that, all
together, lie in some field $F$ which is a finite-dimensional vector space
over $\mathbb{Q}$. The algebraic numbers can then themselves be understood as
rational linear transformations of $F$. Therefore, possibly by passing to
larger
matrices, we can assign to each matrix $\alpha(g_i)$ a larger matrix with
rational entries, and this assignment extends to a monomorphism $\beta$.
Thus, $\beta\alpha(h)$ is a rational matrix assigned to $h$ for every $h\in
G$. 

To define the third and last homomorphism, let $p$ be a prime which
does not divide the denominator of any $\beta\alpha(g_i)$. (The prime
$p$ therefore also does not divide the denominator of any coefficient
of any $\beta\alpha(h)$.) In general, if $p$ does not divide $b$, the
fraction $a/b$ is well-defined as an element of $\mathbb{Z}/p$. Therefore
we can reduce all $\beta\alpha(g_i)$'s mod~$p$ to obtain a homomorphism
$\gamma$ of $\beta\alpha(G)$ to a group of matrices over $\mathbb{Z}/p$
if $p$ fails to divide all denominators in all $\beta\alpha(g_i)'s$.  We
know that $\beta\alpha(g)$ for the originally-chosen $g$ is not the
identity matrix. Therefore if $p$ is larger than all numerators in the
matrix $\beta\alpha(g)$ as well, then $\gamma\beta\alpha(g)$ will also
be distinct from the identity.
\end{proof}

\begin{corollary}
Every hyperbolic lattice is residually finite. \label{cresidual}
\end{corollary}
\begin{proof}For a $d$-dimensional hyperbolic lattice $G$, the base space
$B=\mathbb{H}^d/G$ is a smooth and closed, hence triangulable, manifold.
Therefore the fundamental group $\pi_1(B)$ is finitely generated. In effect,
$G$ is isomorphic to a finitely generated group of matrices.
\end{proof}

\begin{proof}[Proof of theorem] The proof consists of two parts. In the first
part we show the existence of a $d$-dimensional hyperbolic lattice, and in the
second part, given a number $c$, we show that every hyperbolic lattice contains
a sublattice $G'$ whose girth is greater than  $c$. For the first part, we
quote directly from the introduction to \cite[p.~111]{Bo}:

\begin{quotation} A Clifford-Klein form of a connected and simply
connected Riemannian manifold $M$ is a Riemannian manifold $M'$ whose
universal Riemannian covering (universal covering endowed with the metric
lifted from $M'$) is isomorphic to $M$. The main purpose of this note is to
prove the following:

\begin{theorema} A simply connected Riemannian symmetric space
$M$ always has a compact Clifford-Klein form. Any such form $M'$ has a
finite Galois covering which is proper, unless $M'$ is isomorphic to
$M$.
\end{theorema}

We recall that a Riemannian manifold $X$ is symmetric, in the sense of Cartan,
if it is connected and if every point $x$ in $X$ is an isolated fixed point
of an involutive isometry~$s_x$.
\end{quotation}

Obviously, hyperbolic space $\mathbb{H}^d$ is a simply-connected Riemannian
symmetric space, and a compact Clifford-Klein form of $\mathbb{H}^d$ produces
immediately a lattice of isometries of $\mathbb{H}^d$ as the group of covering
transformations.

For the second part of the proof, assume that $G$ is a lattice of
isometries of $\mathbb{H}^d$ and let $c > 0$. Let $B$ be the quotient
manifold $\mathbb{H}^d/G$ and let $p:\mathbb{H}^d\to B$ be the covering
map. The girth of $G$ is determined by the shortest non-contractible
(unbased) loop in $B$.

If $\lambda$ is a non-contractible unbased loop, $\lambda$ represents a
conjugacy class of elements of $\pi_1(M)$. If $g$ is an element of this
conjugacy class and $N$ is a normal subgroup of $\pi_1(M)$ that does not
contain $g$, then $N$ does not contain any conjugate of $g$ either. Since
$B$ is compact, it admits only finitely many homotopy classes with
loops of length at most $c$, say
$\left\{[\lambda_1],[\lambda_2],\ldots,[\lambda_k]\right\}$. Each
$[\lambda_i]$ represents a conjugacy class of an element $g_i\in\pi_1(M)$.
Let $F=\{g_1,g_2,\ldots ,g_k\}$. By Corollary~\ref{cresidual} and Proposition~\ref{pindex},
$\pi_1(M)$ has a finite-index normal subgroup $N$ that does not intersect
$F$. Thus $N$ does not contain any conjugate of any of the $g_i$'s.

Let $q:\widetilde B\to B$ be the covering of $B$ corresponding to $N$. Since
$B$ is compact and $N$ has finite index, $\widetilde B$ is compact as well.
 By the universality of the covering $p:\mathbb{H}^d\to B$, there exists a
covering $\widetilde p:\mathbb{H}^d\to\widetilde B$, hence $\widetilde B$
determines a $d$-dimensional hyperbolic lattice $\widetilde G$. It is clear
that the girth of $\widetilde G$ is greater than $c$, or, in other words, the
length of every non-contractible (unbased) loop in $\widetilde B$ is greater
than $c$, because $q_*:\pi_1(\widetilde B)\to\pi_1(B)$ is a monomorphism,
$q_*\pi_1(\widetilde B)=N$ and $q$ does not increase the length of any loop.
\end{proof}

\section{Dense $n$-reduced coverings}
\label{sdense}

We begin with an example of an arbitrarily dense 2-saturated lattice
covering of $\mathbb{E}^d$ with unit balls ($d\ge 2$). Let $e_1, e_2,\ldots, e_d$
be an orthonormal basis for $\mathbb{E}^d$. Consider the lattice generated by the
vectors $v_i=ae_i$ for $1\le i\le d-1$ and
$$v_d=\left( 1+\sqrt{1-{{a^2}\over 4}(d-1)}\right) e_d +
{a\over 2}\sum_{i=1}^{d-1}e_i,$$
where $\displaystyle 0<a<{2\over\sqrt{d-1}}$.

Clearly, the unit balls centered at the lattice points form a
(simply-)reduced covering. Moreover, each ball covers pairs of points
not contained in any other ball such that the distance between them
approaches 2. Therefore, if two balls are deleted, then one can find
four uncovered points that form the vertices of a parallelogram with
two sides of lengths approaching 2. Since no such set of four points
can be covered by a single unit ball, the covering is 2-reduced.
However, the density of the covering is arbitrarily large for
sufficiently small $a$. Figure~\ref{figure1} illustrates this covering
for $d=2$.

\fig{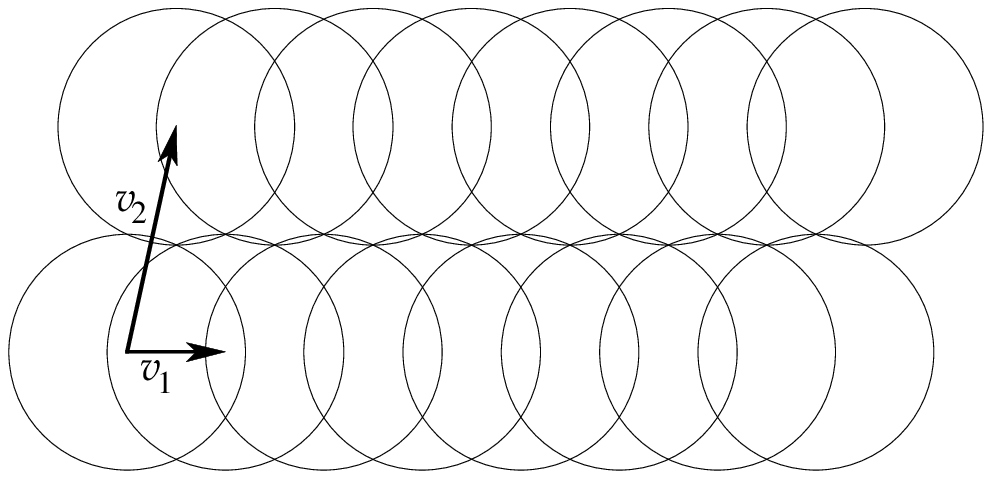}{A high-density, 2-reduced covering by circles.}

An equally simple construction yields an infinitely dense 2-saturated
covering of $\mathbb{E}^d$ with unit balls ($d\ge 2$), in fact locally
infinitely dense at every point. Let $P$ be a hyperplane in $\mathbb{E}^d$
containing the
origin and let $Q$ be a dense subset of $P$ such that $P\smallsetminus Q$ is
dense in $P$ as well. Let $v$ denote the vector normal to $P$ of length 2.
The collection of unit balls centered at the points of the form $q+2iv$ for
$q\in Q$ and $p+(2j+1)v$ for $p\in P\smallsetminus Q$, $i,j\in \mathbb{Z}$, is a
covering. Since the removal of any single ball from this collection uncovers
both end points of the the ball's diameter parallel to $v$, the covering is
not only reduced, but 2-reduced as well.

The second construction generalizes to an arbitrary
body by choosing the hyperplane $P$ to be perpendicular to a diameter of the
body and giving $v$ the same length as the diameter. 

Let $\Theta_n(K)$ be the supremum  of the densities of all $n$-reduced
coverings with $K$. For example, the above constructions show that
$\Theta_1(K)=\Theta_2(K)=\infty$. However, the simple relation
$$\lim\limits_{n\to\infty}\Theta_n(K)=\vartheta(K)$$
implies the existence of a smallest
positive integer $l(K)$ such that $\Theta_n(K)<\infty$ for all $n\ge l(K)$.
The notion of the Newton covering number of a convex body
(introduced below) yields an upper bound for $l(K)$, where $K$ is an
arbitrary convex body in $\mathbb{E}^d$.  In addition, an application of a
theorem of
B\'ar\'any yields a slightly better bound for $l(B^d)$, where $B^d$ denotes
the unit ball in $\mathbb{E}^d$.

Recall that the Hadwiger covering number $H_c(K)$ of a convex body $K$ in
$\mathbb{E}^d$ is the minimum number of translates of $K$ whose union contains a
neighborhood of $K$. Hadwiger \cite{H} asks for the maximum value of $H_c(K)$
over all convex bodies $K$ in $\mathbb{E}^d$. The problem was stated
independently by others, also in the context of the equivalent problem of
illumination of $K$, and it is conjectured that $H_c(K)\le 2^d$ with equality
for parallelotopes only. The conjecture is still open in every dimension
$d\ge 3$.

Similarly, one can consider $N_c(K)$, the minimum number of replicas of $K$
whose union
contains a neighborhood of $K$. The quantities $N_c(K)$ and $H_c(K)$ have
their dual counterparts $N(K)$ and $H(K)$, called the Newton (or kissing)
number and the Hadwiger number, respectively, in the context of packings (see
\cite{FTK}). We call $N_c(K)$ the {\it Newton covering number} of $K$ to
extend the analogy. Obviously, $N_c(K)\le H_c(K)$ for all $K$.

The following theorem establishes a relation between $l(K)$ and $N_c(K)$.

\begin{theorem} For every convex body $K$ in Euclidean
space,
$$l(K)\le N_c(K)+1.$$ \label{thnewton}
\end{theorem}

In the proof of this theorem, as well as in the next section, we use certain
relations between the (global) density of an arrangement and its density with
respect to some bounded domain. In what follows, the volume of a (measurable)
set $S$ will be denoted by $V(S)$. As usual, $B^d$ denotes the unit ball in
$\mathbb{E}^d$, so $rB^d$ is the ball of radius $r$ centered at the origin. Let
$\mathcal A$ be a locally finite arrangement of uniformly bounded measurable
sets, and let $G$ be a bounded domain. The density $d\left(\mathcal A|G\right)$
of $\mathcal A$ relative to $G$ is defined by
$$d\left({\mathcal A}|G\right)={1\over{V(G)}}\sum_{A\in\mathcal A}V(A\cap G),$$
and the average density $d_{av}\left(\mathcal A|G\right)$ of $\mathcal A$ relative to
all translates of $G$ is defined by
$$d_{av}\left({\mathcal A}|G\right)= {\lim_{r\to\infty}}
{1\over{V\left(rB^d\right)}}\int_{rB^d}d\left({\mathcal A}|(G+x)\right),$$
provided the limit exists. Otherwise we take the lim sup or the lim inf,
and it is usually clear from context which limit is meant. Also, the domain of
the integral above is $rB^d$, just as $rB^d$ is frequently used to define
the density of ${\mathcal A}$ as a limit.

The following proposition is derived from these definitions by routine methods
of real analysis, interchanging sums and limits with integrals, and applying
Fubini's theorem.

\begin{proposition} For every locally finite arrangement $\mathcal A$ of uniformly
bounded measurable sets, and any bounded domain $G$, the average density
$d_{av}\left({\mathcal A}|G\right)$ coincides with the density of $\mathcal A$.
\label{paverage}
\end{proposition}

As a direct corollary, we get:

\begin{proposition} Let $\mathcal A$ be a locally finite arrangement of uniformly
bounded measurable sets, and let $G$ be a bounded domain. If the density of
$\mathcal A$ is $a$, then there exists a translate of $G$ such that
$d\left({\mathcal A}|G\right)\ge a$ and there exists a translate of $G$ such that
$d\left({\mathcal A}|G\right)\le a$.
\label{pexists}
\end{proposition}

Proposition~\ref{paverage} can be put in an equivalent, discrete form:

\begin{proposition} Let $D$ be a locally finite set of points, and let $G$ be
a bounded domain. If the number density of $D$ is $a$, then the average
number of points contained in a translate of $G$ is equal to $aV(G)$.
\label{pdaverage}
\end{proposition}

Again, as a corollary, we get:

\begin{proposition} Let $D$ be a locally finite set of points, and let $G$ be
a bounded domain. If the number density of $D$ is $a$, then there exists a
translate of $G$ which contains at least $aV(G)$ points of $D$ and there
exists a translate of $G$ which contains at most $aV(G)$ points of $D$.
\label{pdexists}
\end{proposition}
\begin{proof}[Proof of theorem] By definition, there exists an $\varepsilon>0$
such that some $N_c(K)$ replicas of $K$ cover the $\varepsilon$-neighborhood
(the outer parallel domain of radius~$\varepsilon$) of $K$. Let $p \in K$, and
for each member $K_i$ of the covering, let $p_i$ be the image of $p$ under an
isometry that takes $K$ to $K_i$.   Since the group of isometries of space that
fix $p$ is compact, it can be partitioned into a finite collection of sets such
that, if $g$ and $h$ belong to the same set, then the Hausdorff distance
between $g(K)$ and $h(K)$ is smaller than $\varepsilon /2$, {\it i.e.,} each of
$g(K)$ and $h(K)$ lies in the other's $\varepsilon /2$-neighborhood.

If the covering is sufficiently dense, then by Proposition~\ref{pdexists}, there
exists a ball of radius $\varepsilon /2$ containing $N_c(K)+1$ points $p_i$
such that the Hausdorff distance between each two of the sets $K_i-p_i$ (each
set $K_i$ shifted so that $p_i$ is moved back to the origin) is smaller than
$\varepsilon /2$. We have now $N_c(K)+1$ replicas of $K$, say $K_1, K_2,
\ldots , K_{N_c(K)+1}$, such that $K_i$ lies in the
$\varepsilon$-neighborhood of $K_1$ for $2\le i\le N_c(K)+1$. By the
definition of the Newton covering number of $K$, these $N_c(K)+1$ replicas can
be replaced by $N_c(K)$ others without uncovering any points.
Thus, every sufficiently dense covering with replicas of $K$ fails to
be $\left( N_c(K)+1\right)$-reduced.
\end{proof}

Since $N_c(B^d)=d+1$, the above theorem implies immediately that $l(B^d)\le
d+2$. However, using a result of B\'ar\'any \cite[Th. 2]{Ba} which
generalizes a theorem of Erd\H os and Szekeres \cite{ES}, one can improve
this inequality as follows.

\begin{theorem} $l(B^d)\le d+1$.  \label{thbarany}
\end{theorem}
\begin{proof} B\'ar\'any's theorem states: For any $\varepsilon>0$ and
$d\ge 2$ there exists a constant $n(d,\varepsilon)$ such that every finite
set $V\subset \mathbb{E}^d$ contains a subset $W\subset V$, $|W|\le
n(d,\varepsilon)$ with the property that for $v\in V\smallsetminus W$ there
are points $w_1, w_2 \in W$ with $\measuredangle(w_1vw_2)>\pi -\varepsilon$.
Given positive numbers $\varepsilon$ and $\delta<1$, if a covering of
$\mathbb{E}^d$ with unit balls is sufficiently dense, then, by
Proposition~\ref{pdexists}, some
ball of radius $\delta$ contains at least ${(d+1){n(d,\varepsilon)\choose
2}}$ centers of the unit balls. If we let $V$ be the set of these centers and
we apply B\'ar\'any's theorem, we obtain $d+3$ distinct points $w_1,
w_2,v_1, v_2,\ldots , v_{d+1}$ in $V$ such that
$\measuredangle(w_1v_iw_2)>\pi -\varepsilon$ for $i=1,2,\ldots,d+1$. It
follows that the points $v_1, v_2,\ldots , v_{d+1}$ lie in the ``double
cone'' $C$ (the union of two congruent non-overlapping right cones with a
common base) whose apexes are $w_1$ and $w_2$ and whose angle at each apex is
$2\varepsilon$.

Let $B(p)$ be the unit ball centered at the point $p$. Each of the $d+1$
balls $B(v_i)$ $(1\le i\le d+1)$ is contained in the outer parallel domain
$P$ of radius 1 of $C$. Observe that the set $P\smallsetminus\left(
B(w_1)\cup B(w_2)\right)$ is a neighborhood of the $d-2$-dimensional unit
sphere centered at the midpoint of $w_1w_2$ which lies in the hyperplane
perpendicular to $w_1w_2$, and that his neighborhood is arbitrarily close
to the
sphere for sufficiently small $\varepsilon$ and $\delta$. Since such a
neighborhood can be covered by $d$ unit balls, it follows that a very
dense covering of $\mathbb{E}^d$ with unit balls cannot be $(d+1)$-reduced.
Figure~\ref{figure2} shows the set $P\smallsetminus\left( B(w_1)\cup B(w_2)\right)$  in
dimension 2, where $\varepsilon=\pi /6$ and $\delta=1$ are small enough for
our purpose.
\end{proof}

\fig{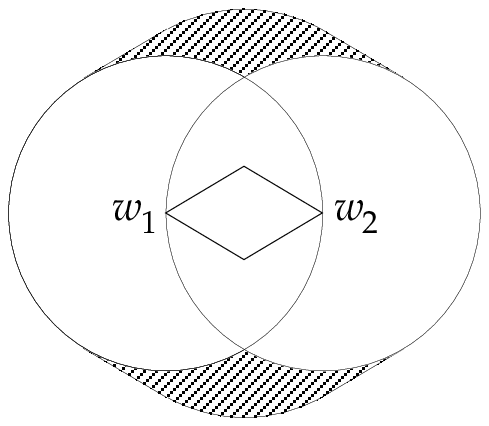}{The set $P\smallsetminus\left( B(w_1)\cup B(w_2)\right)$}

\begin{remark} The proof does not use the fact that the given collection of
balls is a covering. Define an $n$-reduced arrangement (not necessarily a
covering) of replicas of $K$ by the property that it is not possible to
delete $n$ members of the arrangement and replace them with some $n-1$
replicas of $K$ without uncovering any point that was covered by the original
arrangement. Then the above argument demonstrates that a $d+1$-reduced
arrangement of unit balls cannot have arbitrarily high density, and for this
generalization the bound of $d+1$ is the best possible. \end{remark}

So far we have considered arrangements of replicas of a given body, without
restrictions on the isometries that send one of them onto another. However,
we can also consider arrangements with restrictions on the allowed
isometries.  For example, many of the ideas and results of this section
and previous sections generalize to arrangements of translates of $K$
(which is only an appropriate restriction in the Euclidean case, of 
course). In particular, Theorem~\ref{thn} becomes:

\begin{theorem} Every body $K$ in $\mathbb{E}^d$ admits both an $n$-saturated
packing and an $n$-reduced covering in the class of packings and coverings
with translates of $K$.  \label{thtranslates}
\end{theorem}

Theorem~\ref{thnewton} also has a ``translative'' version with a completely analogous
proof. (In the corresponding notation, we indicate the restriction to
translates by the subscript ``$T$''.)

\begin{theorem} For every convex body $K$ in Euclidean
space,
$$l_T(K)\le H_c(K)+1.$$ \label{thnewtranslates}
\end{theorem}

{\sc Example.} While the bounds given in Theorem~\ref{thnewton} and
Theorem~\ref{thbarany} could be far from optimal, the following example
indicates that Theorem~\ref{thnewtranslates} is close to optimal in at least
some cases.

Let $P$ denote a right pyramid in the coordinate space $\mathbb{E}^d$ whose base
is the unit $d-1$-dimensional cube $Q^{d-1}=\{(x_1, x_2,\ldots ,
x_{d-1},0)\in \mathbb{E}^d: |x_i|\le 1/2 {\mbox{ for }} i<d\}$. Observe that
$H_c(P)=1+H_c(Q^{d-1})=1+2^{d-1}$. Consider the set $\Lambda$ of vectors of
the form $(n_1,n_2,\ldots ,n_{d-1},x)$ where each $n_i$ is an integer and
$x\in\mathbb{R}$ is rational if and only if $\sum n_i$ is even. The translates of
$P$ by all vectors of $\Lambda$ is a covering of $\mathbb{E}^d$ of infinite
density. It is easily verified that this covering is $\left(
H_c(P)-1\right)$-reduced. Therefore $l_T(P)\ge H_c(P)$. 

The following estimate for the Hadwiger covering number of a convex body in
$\mathbb{E}^d$ is due to Rogers (unpublished):
$$H_c(K)\le{V(K-K)\over{V(K)}}(d\log d + d\log\log d +5d),$$
where $K-K$ is the difference body of $K$, consisting of points of the form
$x-y$ where $x,y\in K$. The inequality follows from the result of Rogers
\cite{R1} which states that each $d$-dimensional convex body $K$ admits a
covering of $\mathbb{E}^d$ by its translates of density at most $d\log d +
d\log\log d+5d$. If $\{K+a_i\}$ is a covering of density guaranteed
by the theorem of Rogers, then Proposition~\ref{pdexists} implies that there exists a
$\lambda >1$ and a translate $(\lambda K)-K+c$ of $(\lambda K)-K$ containing
at most $(d\log d + d\log\log d+5d)V(K-K)/V(K)$ of the points $a_i$.
Since $(K+a_i)\cap \lambda K+c\ne \varnothing$ if and only if
$a_i\in\lambda K-K+c$, it follows that the respective translates of $K$ cover
$\lambda K+c$.

For a centrally symmetric body $K$ in $\mathbb{E}^d$, we have $V(K-K)/V(K)=2^d$,
thus in this case, the Rogers bound for $H_c(K)$ is reasonably close to
the conjectured best upper bound of $2^d$. In the general case, a result
of Rogers and Shephard \cite{RS} states that $V(K-K)/V(K)\le{{2d}\choose d}$
for every convex body $K\subset \mathbb{E}^d$, which yields the asymptotic bound
$H_c(K)\le 4^{d+o(d)}$. 

\section{Asymptotic density bounds}
\label{sasymptotic}

In Section~\ref{sdense}, we defined $\Theta_n(K)$ as the supremal density of
all $n$-reduced coverings with replicas of $K$, and we mentioned
the simple relation
$$\lim\limits_{n\to\infty}\Theta_n(K)= \vartheta(K).$$
Analogously, let $\Delta_n(K)$ be the infimum of the densities of all
$n$-saturated packings with replicas of $K$, and note the analogous
relation
$$\lim\limits_{n\to\infty}\Delta_n(K)=\delta(K).$$  Also, observe that
$\Delta_n(K)>0$ for every body $K$ and every $n\ge 1$. Obviously, each of the
two sequences $\left\{\Delta_n(K)\right\}$ and $\left\{\Theta_n(K)\right\}$
is monotonic. The following inequalities give estimates for the rate of
convergence of the sequences $\left\{\Delta_n(K)\right\}$ and
$\left\{\Theta_n(K)\right\}$:
\renewcommand{\theequation}{\thesection .\arabic{equation}}
\begin{equation} \Delta_n(K)\ge \delta(K)-O(n^{-1/d}) \label{tweedledum}
\end{equation}
and
\begin{equation} \Theta_n(K)\le \vartheta(K)+O(n^{-1/d}). \label{tweedledee}
\end{equation}

To prove inequality (\ref{tweedledum}), assume that $K$ is a body of
diameter 1 and volume $V$ and let $r$ denote the minimum radius of a ball
that can intersect $n$ non-overlapping replicas of $K$. Let $\sigma_d$ denote
the volume of the unit ball in $\mathbb{E}^d$. By the definition of the
packing density of $K$ and Proposition~\ref{pexists},
$$nV/(\sigma_d r^d)\ge \delta(K)-\varepsilon$$
for every
$\varepsilon>0$, hence $${{nV}\over{\sigma_d r^d}}\ge\delta(K).$$

Assume now that $\mathcal P$ is an $n$-saturated packing with replicas of $K$, and
let $p$ denote the density of this packing. Any ball of radius
$r+2$ must contain at least $n$ members of $\mathcal P$, for otherwise the members
of
$\mathcal P$ contained in the ball could be replaced by $n$ non-overlapping
replicas of $K$ intersecting the concentric ball of radius $r$. Thus the
total volume of the intersections of such a ball with all members of $\mathcal P$
is at least $nV$. Using Proposition~\ref{pexists} again, we obtain
$$p\ge{{nV}\over{\sigma_d (r+2)^d}}.$$
It follows immediately that
$$p\ge\delta(K)\left({{r}\over{r+2}}\right)^d.$$

By the definition of $r$, a ball of radius $r$ intersects at least $n$
non-overlapping replicas of $K$. Each of these replicas is contained in the
concentric ball of radius $r+1$.  Thus $\sigma_d(r+1)^d\ge nV$, and we get
$$r\ge\left({V\over\sigma_d}n\right)^{1/d}-1.$$
Since the function $f(x)=(x/(x+2))^d$ is increasing, we get:
$$p\ge\delta(K)\left({{cn^{1/d}-1}\over{cn^{1/d}+1}}\right)^d,$$
where $c=(V/\sigma_d)^{1/d}$, and inequality (\ref{tweedledum}) follows.

The proof of inequality (\ref{tweedledee}) is analogous.

The above method can be refined as follows to yield some specific density
bounds for $n$-saturated packings and $n$-reduced coverings of $\mathbb{E}^d$
with unit balls. For packings, consider a ``cluster'' of $n$
non-overlapping unit balls and let $G$ be the outer parallel domain of radius
1 of their union. Let $\mathcal P$ be an $n$-saturated packing with unit balls.
Then every translate of $G$ contains at least $n$ centers of the members of
$\mathcal P$. It follows, by Proposition~\ref{pdexists}, that the density of $\mathcal P$ is
at least $n\sigma_d/V(G)$.

The smaller the volume of $G$, the greater the resulting bound, which raises
the problem of arranging $n$ non-overlapping unit balls in $\mathbb{E}^d$ so that
the volume of the outer parallel domain of radius 1 of their union is
minimum. A similar method can be used for coverings, and it leads to the
problem of arranging $n$ unit balls in $\mathbb{E}^d$, this time allowing
overlaps, so that the volume of the inner parallel domain of radius 1 of
their union is maximum. Of course, this method only works if $n\ge d+1$.

\fig{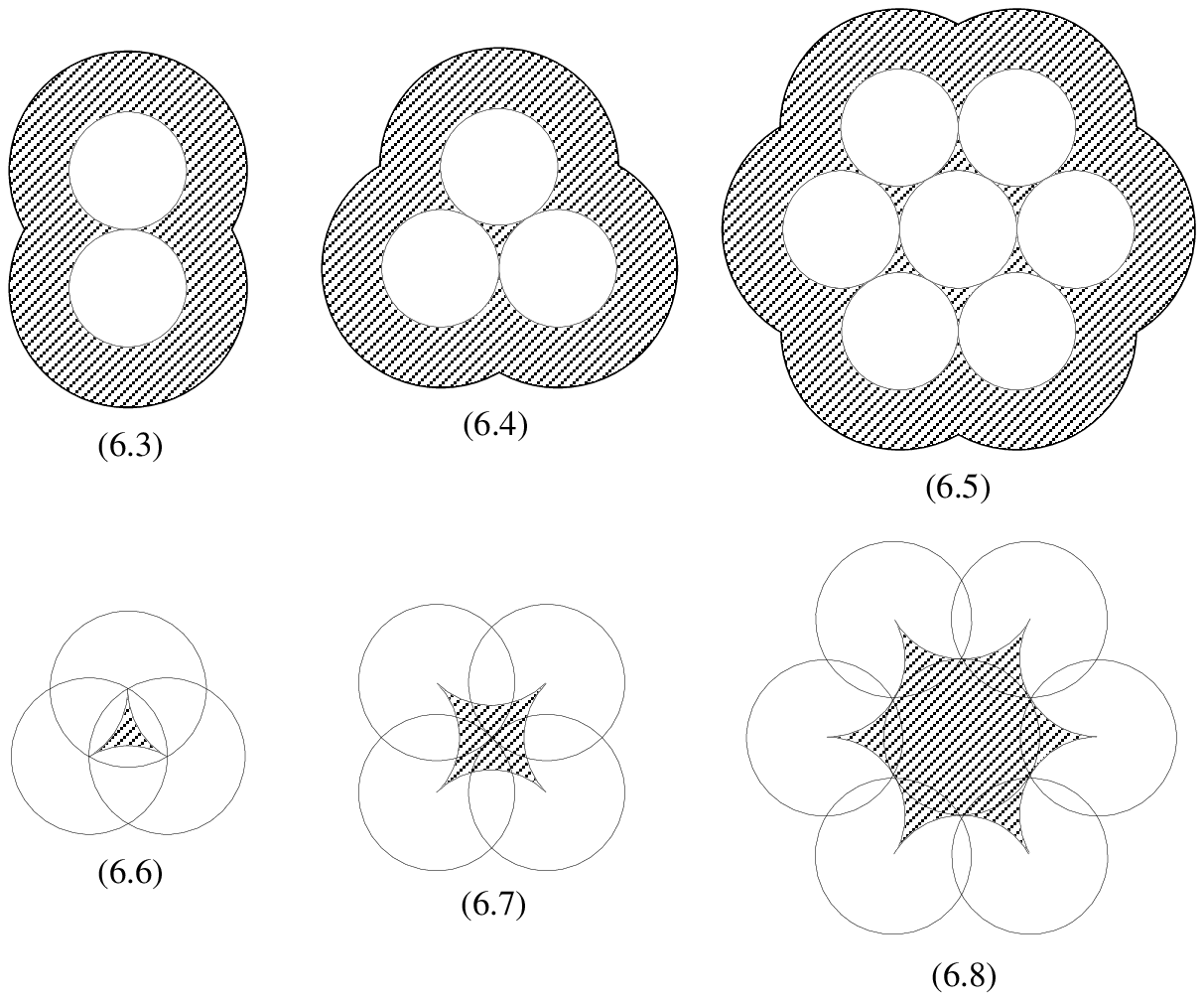}{Economical clusters of unit circles}

Except for some trivial cases, we do not know the solution to these problems
even for $d=2$. However, some clusters of unit circles in $\mathbb{E}^2$, shown
in Figure~\ref{figure3}, seem reasonably economical for the method described above. Using
translates these clusters to estimate
density bounds, we get the following: 
\begin{eqnarray}
\Delta_2(B^2)\ge{{3\pi}\over{3\sqrt 3+8\pi}} & = & 0.31075\dots \\
\Delta_3(B^2)\ge{{3\pi}\over{4\sqrt 3+6\pi}} & = & 0.36561\dots \\
\Delta_7(B^2)\ge{{7\pi}\over{12\sqrt 3+8\pi}} & = & 0.47892\dots \\
\Theta_3(B^2)\le{{6\pi}\over{2\sqrt 3-\pi}} & = & 58.44661\dots \\
\Theta_4(B^2)\le{{4\pi}\over{4-\pi}} & = & 14.63916\dots  \\
\Theta_7(B^2)\le{{7\pi}\over{6\sqrt 3- 2\pi}} & = & 5.35179\dots\ \ .
\end{eqnarray}

Clearly, these inequalities are far from sharp.  Among good estimates
for $\Delta_n(B^2)$ and $\Theta_n(B^2)$, one
stands out. Clearly, any 1-saturated packing with unit balls
becomes a covering if the radius of each ball is increased to 2. Since
$\vartheta(B^2)=2\pi/\sqrt{27}$ (a well-known result of Kershner \cite{K}),
it follows that
$$\Delta_1(B^2)=\pi /6\sqrt {3}=0.302299\dots.$$
Also, as we mentioned before,
$$\Theta_1(B^2)=\Theta_2(B^2)=\infty.$$
Apart from these three cases, it seems difficult to determine
the exact values of $\Delta_n(B^2)$ and $\Theta_n(B^2)$.

\section{Remarks, open problems and conjectures}
\label{sopen}

In relation to the conjecture stated in the introduction, claiming the
existence of completely saturated packings and completely reduced coverings,
observe the following:

\begin{description}

\item[(i)] Complete saturation implies maximum density and complete reduction
implies minimum density. More precisely, the density of a completely saturated
packing with
replicas of a body $K$ exists and is equal to $\delta (K)$. Similarly, the
density of a completely reduced covering with replicas of $K$ exists and is
$\vartheta (K)$.

\item[(ii)] Obviously, the converse of (i) is false. But a weaker statement
holds: A periodic packing with replicas of $K$ with density $\delta(K)$
is completely saturated and a periodic covering with replicas of $K$ whose
density is $\vartheta(K)$ is completely reduced.

\end{description}

The first observation indicates that the conjecture on existence of completely
saturated packings and reduced coverings is not as obvious as it might
appear.  The conjecture, if true, would imply a version of Groemer's
result~\cite{G} on the existence of maximum density packings and
minimum density coverings.

The second observation brings to mind the well-known problem: Given a body $K$,
is there a periodic packing [covering] with replicas of $K$, whose density is
$\delta(K)$ [$\vartheta (K)$]? A positive answer to this question would imply
our conjecture. However, Schmitt \cite{S1} constructed a strictly star-shaped
prototile for a monohedral tiling in $\mathbb{E}^3$ such that no tiling with its
replicas is periodic, and by a slight modification of Schmitt's construction
Conway produced a convex prototile with this property. For $\mathbb{E}^2$ no such
example is known, but according to another result of Schmitt \cite{S2} there
is a strictly star-shaped set $K\subset\mathbb{E}^2$ whose replicas do not admit
a periodic packing of density $\delta(K)$.

Generally, it seems extremely difficult to determine whether a given convex
body admits a periodic packing (covering) of maximum (minimum) density. In
particular, the answer is not known for the $d$-dimensional ball ($d\ge 3$).
The case $d=2$ offers some answers, since it is known (see \cite{FT1,FT3})
that every centrally-symmetric convex disk attains its packing density in a
lattice packing. The analogous statement for coverings is only a conjecture,
supported by a partial result under the restriction to crossing-free
coverings (see \cite{FT1,FT3}).

There are only a handful of cases in which sphere packings in Euclidean or
hyperbolic space are known to be completely saturated. Without exception,
they follow from the Rogers \cite{R2} and the B\"or\"oczky \cite{Bor2}
bounds:  The density of any sphere packing in $d$ dimensions is at most the
density in a regular simplex of $d+1$ kissing spheres with centers at the
vertices of the simplex.  If the regular simplex tiles space, there exists a
corresponding periodic sphere packing that achieves the bound.  The only
regular simplices that tile Euclidean and hyperbolic space are:

\begin{description}

\item[$\bullet$] Equilateral triangles in $\mathbb{E}^2$.

\item[$\bullet$] Triangles in $\mathbb{H}^2$ with angles of $2\pi/n$ for $n\ge 7$.

\item[$\bullet$] Simplices in $\mathbb{H}^4$ with dihedral angles of $2\pi/5$.

\end{description}

The analogous bound for coverings also holds in Euclidean space (see
Coxeter-Few-Rogers \cite{CFR}), although it is open for sphere
coverings of hyperbolic $d$-space and the $d$-sphere for $d>2$.
(For $\mathbb{H}^2$ the bound for circle coverings is due to L.~{Fejes
T\'oth} \cite{FT2}.)  Therefore the same simplices also produce
completely reduced coverings of the same types, except perhaps in
$\mathbb{H}^4$.

\fig{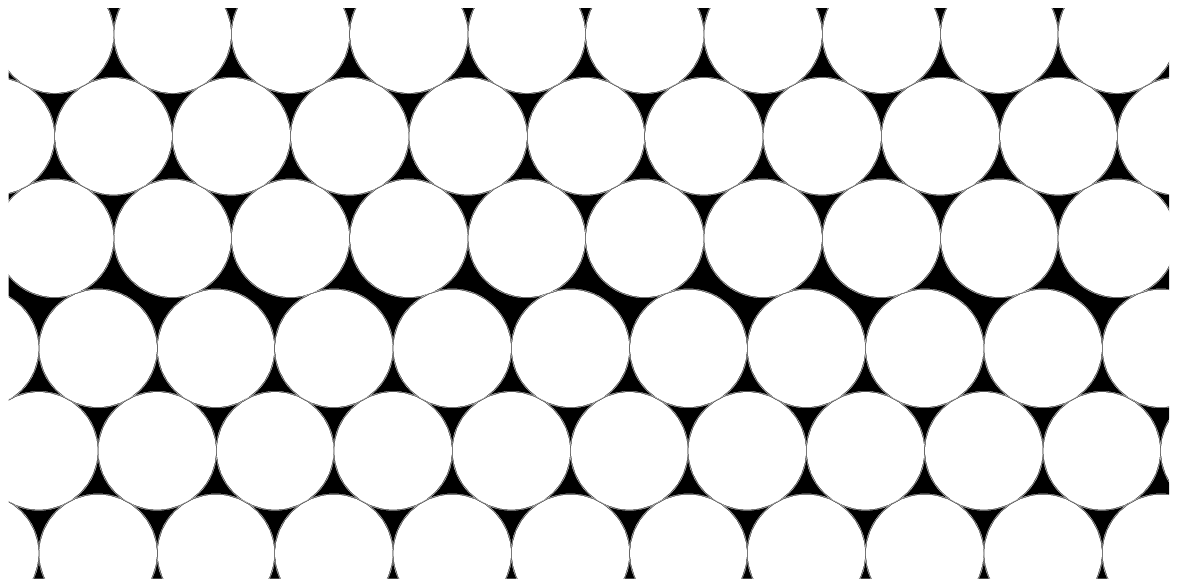}{A sub-optimal packing which could be completely saturated}

Clearly, the familiar densest lattice packing (covering) of $\mathbb{E}^2$
with unit circles is completely saturated (reduced), but we do not know
whether there is a non-lattice, completely saturated packing
(completely reduced covering) with unit circles. We do not even know
whether or not the circle packing in Figure~\ref{figure5} is completely saturated.
The arrangement of circles shown there is given by dividing lattice
packing into two ``half-plane'' parts along a pair of adjacent rows of
circles and then separating the parts slightly while maintaining
contact between the two adjacent rows.

Although it seems difficult to determine $\Delta_n(B^2)$ and $\Theta_n(B^2)$,
we conjecture that
$$\Delta_2(B^2) = \pi(3-\sqrt5)/\sqrt{27}=0.461873\dots,$$
which is the density of the packing shown in
Figure~\ref{figure6}. But we do not even have a conjecture for the other constants.

\fig{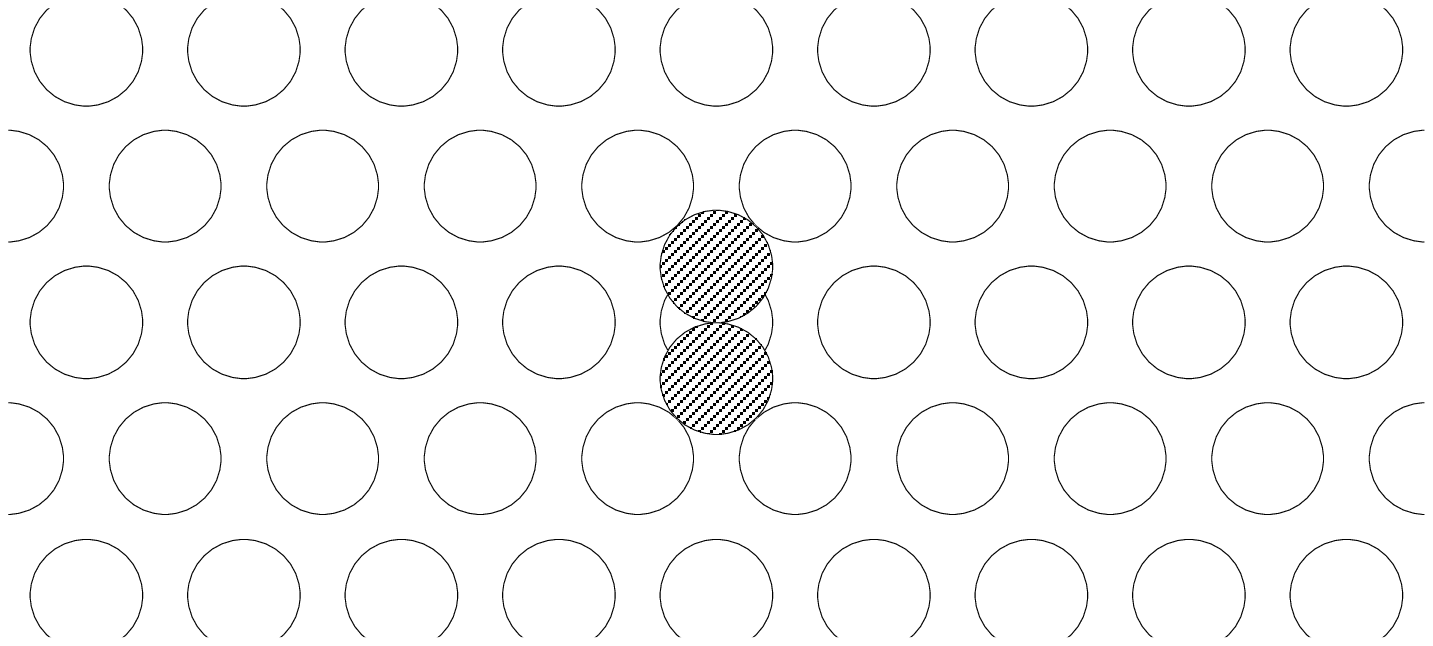}{Possibly the least dense 2-saturated packing}

Theorem~\ref{thnewton} relates the Newton covering number $N_c(K)$ for a
convex body $K$ to $l(K)$. Newton covering numbers are of interest in their
own right: Among all
convex bodies in $d$ dimensions, which one has the greatest Newton covering
number, and what is that number? Let $P$ denote the right pyramid over a
$d-1$-dimensional cube $Q^{d-1}$, as in the remark following
Theorem~\ref{thnewtranslates}. Is
$N_c(P)=N_c(Q^{d-1})+1$? Does $N_c(P)$ depend on the height of the pyramid?
What is the Newton covering number of the cube $Q^d$?

The inequality $l(B^d)\le d+1$ (Theorem~\ref{thbarany}) is sharp for $d=2$
(see Figure~\ref{figure1}), but we suspect that for $d\ge 3$ this is not the case. It
might even turn out that $l(B^d)=3$ for all $d\ge 2$. This problem can be
stated more simply as follows: Given a very dense covering of $\mathbb{E}^d$
with unit balls, can one always make a new covering by replacing three balls
by two?

\begin{flushleft}
Gabor Fejes T\'oth: {\sc Mathematical Institute of the Hungarian
Academy of Sciences, P.O.Box 127, Budapest, H-1364, Hungary.} \\
{\it E-mail address:} gfejes@math-inst.hu \\
\mbox{ } \\

Greg Kuperberg:  {\sc Department of Mathematics, Yale University,
New Haven, CT 06520, U.S.A.} \\
{\it E-mail address:} greg@math.yale.edu \\
\mbox{ } \\

W{\l}odzimierz Kuperberg:  {\sc Department of Mathematics, Auburn University,
AL~36849-5310, U.S.A.} \\
{\it E-mail address:} kuperwl@mail.auburn.edu
\end{flushleft}

\end{document}